\documentclass{article}
\usepackage{amsfonts,amssymb}
\newcommand{\C}{\mathbb{C}}
\newcommand{\N}{\mathbb{N}}
\newcommand{\R}{\mathbb{R}}

\newcommand{\B}{\mathfrak{B}}
\newcommand{\A}{\mathfrak{A}}

\newcommand{\fC}{\mathfrak{C}}
\newcommand{\fN}{\mathfrak{N}}

\newcommand{\cL}{{\cal{L}}}
\newcommand{\cM}{{\cal{M}}}
\newcommand{\cS}{{\cal{S}}}

\date{}

\newtheorem{tht}{Theorem}[section]
\newtheorem{thd}[tht]{Definition}

\newtheorem{thp}[tht]{Proposition}

\newtheorem{thex}[tht]{Example}

\title{On representations of partial $^\ast$-algebras based on $\B$-weights}
\author{Klaus-Detlef K\"ursten  and Elmar Wagner
        \thanks{Universit\"at Leipzig, Mathematisches Institut} }
\begin{document}
\maketitle
\begin{abstract}
A generalization of the $GNS$-representation is investigated that represents
partial $^\ast$-algebras as systems of
operators acting on a partial inner product
space ($PIP$-space). It is based on possibly indefinite $\B$-weights which
are closely related to the positive $\B$-weights introduced by
J.-P. Antoine, Y. Soulet and C. Trapani.
Some additional assumptions had to be made in order to  guarantee the
$GNS$-construction.
Different partial products of operators on a $PIP$-space are considered
which allow the $GNS$-construction under suitable conditions. Several
examples illustrate the argumentation and indicate inherent problems.
\end{abstract}
%
%
\section{Introduction}
%
%
The development of the theory of partial $^\ast$-algebras has been
motivated by the appearance of such structures in models of local
quantum field theory and quantum statistical mechanics (e.g., see
\cite{Bor72, Las79}). Lately, several standard results of the theory
of $^\ast$-algebras of operators are extended to a certain degree to 
partial $^\ast$-algebras, for instance representation theory,
modular theory of Tomita-Takesaki, and automorphism groups and
$^\ast$-derivations. For details and further references, we refer 
to the review \cite{AIT96} by J.-P. Antoine, A. Inoue and C. Trapani.

As the $GNS$-construction is one of the basic tools of the theory
of $^\ast$-algebras,
there arises a particular
interest in extending it to
partial $^\ast$-algebras. A promising approach to this problem
has been made by J.-P. Antoine, A. Inoue and C. Trapani \cite{AIT91},
starting with a positive sesquilinear form on a partial $^\ast$-algebra,
using a subspace of the space of all right multipliers to set up
the representation, and taking into account the possible lack of
(semi-)associativity. The result is a representation of the
partial $^\ast$-algebra into the partial 
$O^\ast$-algebra $L^{\dag}(D,H)$.
Nevertheless, this approach might be not general enough.
Yet for $^\ast$-algebras there exists a $GNS$-construction based
on weights, that is, positive functionals defined on the positive cone
of the  $^\ast$-algebra that do not necessarily take finite values.
In order to give a $GNS$-construction for partial $^\ast$-algebras that 
generalizes also the theory of weights,
J.-P. Antoine, Y. Soulet and C. Trapani \cite{AST95}
introduced the notion of a (positive)
$\B$-weight. The $GNS$-construction based on a $\B$-weight will lead
to a representation of the given partial $^\ast$-algebra as a
system of operators 
acting on some partial inner product space ($PIP$-space); 
the basic ideas of $PIP$-spaces were 
developed by J.-P. Antoine and A. Grossmann in earlier 
papers \cite{Ant76a}.

Our treatment of the subject follows closely the program presented in
\cite{AST95},
except that we do not require the $\B$-weight to be positive.
However, there will be included several examples which show that
our results are also relevant to the case of positive $\B$-weights.
In Section 2 we give a sufficient and necessary condition
that a $\B$-weight determines the structure of a non-degenerate
$PIP$-space. Sections 3 and 4 are devoted to the study of a
generalized $GNS$-representation of a partial $^\ast$-algebra
as systems of
operators acting on a non-degenerate $PIP$-space $V$.
To do this, the linear space of all continuous linear operators on $V$
(denoted by $Op(V)$) must be equipped with a multiplicative
structure. In Section 3 the underlying multiplicative structure is
based on a definition due to J.-P. Antoine
and A. Grossmann \cite{Ant76a}. It turns out that additional
assumptions must be made in order to guarantee the
$GNS$-representation. Under certain conditions some of the
additional assumptions can be removed by introducing more general
partial products on $Op(V)$. This is the central theme of
Sec\-\mbox{tion 4}.  For the convenience of the
reader, the necessary definitions concerning $PIP$-spaces,
partial $^\ast$-algebras and $\B$-weights are included in
Sections 2 and 3. For a more detailed study we refer to
the references \cite{Ant76a, AIT96, AIT91, Ant85, AST95, Kue90b}.
%
%
\section{Construction of non-degenerate $PIP$-spaces}
%
%
%
In this section we investigate the problem whether, given a (possibly
indefinite) $\B$-weight $\Omega$ on a partial $^\ast$-algebra,
there exists a non-degenerate $PIP$-space
associated  with $\Omega$ in a natural way. It turns out that
the existence of such a $PIP$-space may be characterized by  additional
conditions on the $\B$-weight. Examples show that these conditions do not
follow from the axioms of $\B$-weights.

Let us summarize some notations and definitions concerning $PIP$-spaces,
partial $^\ast$-algebras and $\B$-weights. These
definitions are essentially equivalent to the original definitions
in \cite{Ant76a, AST95},
except that our forms are linear in the first argument and that
the form $\Omega$ is not required to be positive semi-definite here.
A weak linear compatibility on a $\C$-vector space $V$ is a symmetric binary
relation $\#$ on $V$ such that all non-empty sets of the type
$$M^\#\  \stackrel{\rm def}{=}\  \{ \varphi \in V\,;\, \varphi \, \# \,
\psi \ \,\mbox{ for all }
\psi \in M \} \qquad (M \subset V)$$
(called assaying subspaces) are linear subspaces of $V$. A linear
compatibility on $V$ is a weak linear compatibility such that  all
sets of the type
$M^\#$ are non-empty. Suppose now that $\#$ is a linear
compatibility on the $\C$-vector space $V$ and let
$\Gamma(\#)$ denote the graph of $\#$.
A partial inner product on $(V, \#)$ is a mapping
$$\Gamma(\#) \ni (\varphi, \psi) \rightarrow \langle \varphi ,
\psi \rangle \in \C $$
which is linear in $\varphi$ and satisfies $\langle \varphi ,
\psi \rangle = \overline{\langle \psi , \varphi
 \rangle }$  whenever $ \varphi \, \# \, \psi $. A triple
$(V, \#, \langle . ,\! . \rangle)$, where $\langle .,\! .\rangle$ is
a partial inner product on $(V, \#)$, is called a $PIP$-space.
It is said to be non-degenerate if and only if $\langle \varphi ,
\psi \rangle = 0$ for all $\psi \in V^\#$ implies that
$\varphi=0$.

A partial $^\ast$-algebra is a $^\ast$-vector space $\A$ together with
a subset $\Gamma \subset \A \times \A$ and a mapping
$\Gamma\, \ni\, (x, y)\, \rightarrow\, x \cdot y = x y\, \in\, \A$ such that
$(x, y), \, (x, z) \in \Gamma $ and $\lambda, \mu \in \C$ imply that
$(x, \lambda y + \mu z)\in \Gamma$, $(y^\ast, x^\ast)\in \Gamma$, and that
\begin{eqnarray*}
x \cdot (\lambda y + \mu z)&=& \lambda (x \cdot y) + \mu (x \cdot z)\ , \\
(x \cdot y)^\ast &=& y^\ast \cdot x^\ast\ .
\end{eqnarray*}
The set of right multipliers of a subset $\fN \subset \A$ is the set
$$R({\fN})\ \stackrel{\rm def}{=}\  \{ x\in\A\,;\,(y,x)\in\Gamma\ \,
\mbox{ for all } y
\in {\fN} \}\ .$$

$GNS$-constructions map abstract partial $^\ast$-algebras into
certain spaces of linear operators acting on linear spaces, for instance 
on unitary spaces or on $PIP$-spaces.
Such spaces of linear operators represent also basic
examples of spaces with partially defined products. Examples in
\cite{Kue86b, Kue90b, Kue94} show that these spaces are not necessarily
associative. For this reason, associativity is not included in the
axioms of partial  $^\ast$-algebras.

%
%
Now (possibly indefinite) $\B$-weights are defined as follows.
\begin{thd}           \label{Bweight}
Suppose $\A$ is a partial $^\ast$-algebra,
$\B$ is a linear subspace of $R(\A)$, and
$\sharp$ is a weak linear compatibility on $\A$ with graph $\Gamma(\sharp)$.
A mapping
$$\Omega\;:\;\Gamma(\sharp)\,\ni\,(x, y)\,\rightarrow\,\Omega(x,y)\,\in\,\C$$
which is linear in $x$ and satisfies $$\Omega(x,y)\ =\ \overline{\Omega(y,x)}$$
whenever $(x,y) \in \Gamma(\sharp)$
 is said to be a $\B$-weight if the following conditions are satisfied: \\
i) $\B \times \B \cup \A  \B \times \B \subset \Gamma(\sharp) \quad
        (\mbox{where } \A \B \stackrel{\rm def}{=} \{ a b \,;\, a \in \A
        \mbox{ and } b \in \B \})$. \\
ii) $\Omega(xb_1,b_2) =\Omega(b_1,x^\ast b_2) \mbox{ for all }
        x \in \A \mbox{ and } b_1, b_2 \in \B  $. \\
iii) If $ x_1, x_2 \in \A$ and $x_1 \in R(\{x_2\})$, then
     $(x_1 b_1 , x_2^\ast b_2)\in\Gamma(\sharp)$ for all
     $b_1, b_2 \in \B$ and
     $\Omega(x_1 b_1 , x_2^\ast b_2) =\Omega((x_2 x_1) b_1 , b_2)$.\\
iv) If $x \in \B$ and $\Omega(y,x)=0$ for all $y \in \B$, then
        $\Omega(y,x)=0$ for all $y \in \B^\sharp$.
\end{thd}

If $\A$ is a $^\ast$-algebra and $\omega$ a linear functional on $\A$
satisfying $\omega(x^\ast x)\in\R$ for all $x\in\A$, then one obtains
a sesquilinear form $\Omega$ on $\A$ by setting
$\Omega(x,y)=\omega(y^\ast x)$. It arises the question
if we can start with a linear functional on a partial $^\ast$-algebra $\A$
and construct a $\B$-weight in a similar way. Naturally,
additional assumption must be made. In view of
Definition \ref{Bweight} ii) and iii), we impose on $\A$
semi-associativity: $\A$ is called semi-associative if $x\in R(\{y\})$
implies $xb\in R(\{y\})$ for all $b\in R(\A)$ and $y(xb)=(yx)b$.
The following proposition describes a situation where a linear functional
determines the structure of a $\B$-weight.
\begin{thp}    \label{omega}
Suppose $\A$ is a semi-associative partial $^\ast$-algebra, $\B$ is a
linear subspace of $R(\A)$, $D$ is a linear subspace of $\A$,
and $\omega:D\rightarrow \C$ is a linear functional such that
the following conditions are satisfied: \\
i) $D=D^\ast\, ({\stackrel{\rm def}{=}}\ \{x^\ast \,;\, x \in D \})$
and $\omega(x^\ast)=\overline{\omega(x)}$ for all
   $x\in D$, \\
ii) $\B^\ast\B\subset D$ and $\B^\ast\A\B(=\B^\ast(\A\B))\subset D$. \\
iii) If $x\in\B$ and $\B^\ast x\subset ker(\omega)$, then
     ${\B'}^\ast x\subset ker(\omega)$ where $\B'\stackrel{\rm def}{=}
\{$\mbox{$y\in\A$}\,; \mbox{$\B^\ast y\subset D\}$}.\\
Define $x\,\sharp\,y$ if and only if $x\in R(\{y^\ast\})$ and
$y^\ast x\in D$. Then $\sharp$ is a linear compatibility on $\A$,
and $\Omega$, defined by $\Omega(x,y)=\omega(y^\ast x)$ whenever
$x\,\sharp\,y$, is a $\B$-weight on $\A$ in the sense of
Definition \ref{Bweight}.
\end{thp}
{\sc Proof:}
Since $(y^\ast,x)\in\Gamma$ and $y^\ast x \in D$ imply that
$(x^\ast,y)\in\Gamma$ and $x^\ast y \in D$,
$\sharp$ is symmetric.
Since the partial product on $\A$ is distributive and
$D$ is a linear subspace, $\sharp$ defines a linear compatibility.
Using the hypothesis, one proves easily that $\Omega$ satisfies
the conditions of Definition \ref{Bweight}. For instance,
if $x_1\in R(\{x_2\})$, then repeated application of semi-associativity
gives for all $b_1,b_2\in \B$
$$D\,\ni\, b_2^\ast((x_1 x_2)b_1)\ =\  b_2^\ast(x_1(x_2 b_1))
                             \ =\  (b_2^\ast x_1)(x_2 b_1)\ ,$$
which implies Definition \ref{Bweight} iii).
To verify Definition \ref{Bweight} iv), observe that
\mbox{$\B^\sharp =\B'$} and apply iii).
\medskip

The following proposition is the main result of this section.
It characterizes $\B$-weights to which there
is associated  a non-degenerate $PIP$-space in a natural
way. $PIP$-spaces obtained in this manner will serve in the following
sections as representation spaces for a generalized $GNS$-construction.
\begin{thp}           \label{PIPspace}
Let subspaces of $\B^\sharp$ be defined by
\begin{eqnarray*}
{\fN}_1\  =\  \{x \in \B^{\sharp \sharp} \,;\, \Omega(x,y)=0 \ \,
        \mbox{ for all } y \in \B^\sharp \}\ , \\
{\fN}_2\  =\  \{x \in \B^\sharp \,;\, \Omega(x,y)=0 \ \,
        \mbox{ for all } y \in \B^{\sharp \sharp} \}\ .
\end{eqnarray*}
Then for a linear subspace $\fN \subset \B^\sharp$ the following conditions
are equivalent: \\
i) On $V\stackrel{\rm def}{=}\B^\sharp / \fN$ there exists a non-degenerate
$PIP$-space
structure such that $$(x+\fN) \; \# \; (y+\fN) , \mbox{ if and only if }
 x \; \sharp \; y\ ,$$ and
that $$\langle x+\fN,y+\fN\rangle\ =\ \Omega(x,y)$$ whenever $x\;\sharp\;y$.\\
ii) $\fN = {\fN}_1 = {\fN}_2$.
\end{thp}
{\sc Proof:} ii)$\Rightarrow$i):\
Let $x_1, x_2\in \B^\sharp$ and $n_1, n_2\in\fN$. If $(x_1,x_2)\in
\Gamma(\sharp)$, then we observe that $(x_1+n_1,x_2)\in\Gamma(\sharp)$;
using firstly $n_1\in\B^{\sharp\sharp}$ since $\fN=\fN_1$ and secondly
the linearity of $\sharp$.
Continuing in  this way gives $(x_1+n_1,x_2+n_2)\in\Gamma(\sharp)$.
Hence $(x_1+n_1 , x_2+n_2)\in\Gamma(\sharp)$ for all $n_1,n_2\in\fN$
if and only if $(x_1,x_2)\in\Gamma(\sharp)$; thus $\#$ is well defined.
By  the properties of $\Gamma(\sharp)$, it follows  that  $\#$
defines a linear compatibility on $V$. The above arguments allow us to write
$\Omega(x_1+n_1,x_2+n_2)=\Omega(x_1,x_2)+\Omega(x_1,n_2)+\Omega(n_1,x_2)+
\Omega(n_1,n_2)$ whenever $(x_1,x_2)\in\Gamma(\sharp)$.
Since $\fN=\fN_1$, we get $\Omega(x_1+n_1,x_2+n_2)=\Omega(x_1,x_2)$.
Hence $\langle .,\! .\rangle$ is well defined. That $\langle .,\!.\rangle$
defines a partial inner product follows from the assumed properties
of $\Omega$. By $\fN=\fN_2$, $V$ is non-degenerate.\\
i)$\Rightarrow$ii):\
Clearly, $\fN\subset {\B}^{\sharp \sharp}$ is necessary, otherwise $\#$
would not be well defined. Suppose $n\in \fN$, $y\in {\B}^{\sharp}$
and $\Omega(n,y)\neq 0$. Then for any $x\in {\B}^{\sharp\sharp}$ we get
\mbox{$\Omega(x+n,y)\neq$}\,$\Omega(x,y)$,
so $\langle .,\!.\rangle$ is not well defined;
therefore $\fN\subset {\fN}_1$. Note that $V^\# =
\{x+\fN\,;x\in{\B}^{\sharp\sharp}\}$. As $V$ is non-degenerate, we have
${\fN}_2\subset \fN$. Finally, \mbox{${\fN}_1\subset {\fN}_2$} since
${\B}^{\sharp\sharp}\subset{\B}^{\sharp}$. Thus we have proved
${\fN}_1\subset{\fN}_2\subset\fN\subset{\fN}_1$ which yields ii).
\medskip

The following two examples illustrate that for arbitrary $\B$-weights
the assertions of Proposition \ref{PIPspace} are not necessarily
satisfied. The first example shows that the assertion $\fN_1=\fN_2$
depends strongly on the weak linear compatibility $\sharp$.
The second example is especially designed to show that there is no
natural generalization of the Schwarz inequality for positive
$\B$-weights.
\begin{thex}  \label{exI}
There exist a partial $^\ast$-algebra $\A$, a linear subspace
$\B \subset \A$, and a $\B$-weight $\Omega$ on $\A$
such that ${\fN}_1\neq{\fN}_2$.
\end{thex}
Consider the $^\ast$-algebra $\A\stackrel{\rm def}{=}L_{\infty}([0,2])$
of all measurable, bounded, complex
functions on the interval $[0,2]$. Set $\B=
\chi_{{ }_{[0,1]}}\A$,
where $\chi_{{ }_{[0,1]}}$ denotes the characteristic function of the
interval $[0,1]$. Define $\sharp$ by setting
$\Gamma(\sharp)=(\A\times \B)\cup (\B\times \A)$
and define $\Omega(f,g)=\int_{0}^{2} f{\bar g}dt$ whenever
$(f,g)\in \Gamma(\sharp)$.
Using $\A\B\subset\B$, one easily verifies that $\Omega$ is a $\B$-weight
in the sense of Definition \ref{Bweight}. Now ${\B}^{\sharp}=\A$,
${\B}^{\sharp\sharp}=\B$, and thus
${\fN}_1=\{0\}\neq\chi_{{ }_{[1,2]}}\A={\fN}_2$.
Notice that if we had defined $\Gamma(\sharp)=\A\times\A$,
$\fN_1=\fN_2$ would hold.
\medskip

Incidentally, Example \ref{exI} yields an example of
Proposition \ref{omega}, too;  just set $D=\B$ and
define $\omega(f)=\int_{0}^{2} f dt$ for all $f\in D$.

We call a $\B$-weight $\Omega$ positive if
$\Omega(x,x)\ge 0$ whenever $(x,x)\in \Gamma(\sharp)$.
One might try to obtain a better result as
Proposition \ref{PIPspace} by employing the Schwarz inequality.
If, for instance, $\Omega$ is a positive semi-definite 
sesquilinear form on $\A$,
then it is sufficient to consider the set
${\fN}= \{x \in \B^{\sharp} \,;\, \Omega(x,x)=0 \}$
since the Schwarz inequality implies $\Omega(y,n)=0$
for all $y\in\B^{\sharp}$ and $n\in\fN$. Unfortunately, this approach is
useless for $PIP$-spaces with positive $\B$-weights.
It can, namely, happen that $\Omega(y,n)\neq 0$ although
$n\in \B^{\sharp\sharp}$ and $\Omega(n,n)=0$. We shall present an
explicit example.
\begin{thex}
There exist a partial $^\ast$-algebra $\A$, a linear subspace
$\B \subset \A$, and a positive $\B$-weight $\Omega$ on $\A$
such that the following statements hold: There are $b\in
\B^{\sharp\sharp}$ and $y\in \B^{\sharp}$ such that
$\Omega(b,b)=0$ but $\Omega(y,b)\neq 0$,
$\fN_1,\fN_2\subset\B^{\sharp\sharp}$ and
${\fN}_1\neq{\fN}_2$.
\end{thex}
Let $\B$ be the vector space $\varphi$ of all complex sequences
$(x_n)=(x_n)_{n\in\N}$ with a finite number of nonzero entries.
Consider the complex vector space
$$
\A\ \stackrel{\rm def}{=}\ \C(n)+\C(\mbox{$\frac{1}{n}$})+\varphi\ ,
$$
where $(n)$ and $(\frac{1}{n})$ denote the sequences
$(n)_{n\in\N}$ and $(\frac{1}{n})_{n\in\N}$, respectively.
$\A$ becomes a partial $^\ast$-algebra by restricting the
pointwise multiplication of sequences to
$\Gamma \stackrel{\rm def}{=} (\A\times\B)\cup(\B\times\A)$
and defining an involution by complex conjugation.
 We introduce a linear compatibility
$\sharp$ and a positive $\B$-weight $\Omega$ on $\A$ by setting
\begin{eqnarray*}
\Gamma(\sharp) &=& (\A\times(\C(\mbox{$\frac{1}{n}$})+\B))
                 \cup((\C(\mbox{$\frac{1}{n}$})+\B)\times\A)\ ,\\
\Omega((x_n) , (y_n)) &=&
\lim_{k\to\infty} x_k\overline{y_k}\,
\quad {\rm whenever}\ ((x_n),(y_n))\in \Gamma(\sharp)\ .
\end{eqnarray*}
Notice that $\Omega((x_n),(v_n))=0$ for all $(x_n)\in\A$,
$(v_n)\in\B$ and that the product of any pair
$((x_n),(y_n))\in \Gamma$ lies in $\B$.
Combining these two facts, one proves easily that $\Omega$
is indeed a $\B$-weight. Moreover, $\B^{\sharp}=\A$ and
$\B^{\sharp \sharp}=\C(\frac{1}{n})+\B$. \\
Given $(v_n),(w_n)\in\B$, $\alpha,\beta_1,\beta_2\in\C$, we calculate
$$
\Omega(\beta_1(\mbox{$\frac{1}{n}$})+(v_n),\alpha(n)+
\beta_2(\mbox{$\frac{1}{n}$})+(w_n))\ =\ \beta_1\overline{\alpha}\ .
$$
As an example, the choice $b=(\frac{1}{n})\in \B^{\sharp\sharp}$,
$y=(n)\in \B^{\sharp}$ gives
$\Omega(b,b)=0$ and $\Omega(y,b)=1\neq 0$. Furthermore,
it follows ${\fN}_1=\B$ and ${\fN}_2=\B^{\sharp \sharp}$, hence
${\fN}_1\neq{\fN}_2$ as asserted. Note that
${\fN}_1,{\fN}_2\subset\B^{\sharp\sharp}$ is satisfied.
%
%
%
\section{Representations of partial $^\ast$-algebras
                      using factorization  products}
%
%

We turn now to the problem of constructing a representation of partial
$^\ast$-algebras as systems of 
operators acting on a non-degenerate $PIP$-space.
The non-degenerate $PIP$-space in question will be, of course, the one
constructed in the preceding section. For non-degenerate $PIP$-spaces we can
adopt the concepts of dual pairings; if $(V,\#,\langle.,\!.\rangle)$
is a non-degenerate $PIP$-space and $X\subset V$ an assaying subspace,
then the restriction of $\langle.,\!.\rangle$ to
$X\times \overline{X^\#}$, where $X$ is an assaying subspace and
$\overline{X^\#}$ denotes the
associated conjugate linear space of $X^\#$, is a dual pairing.
From now on, each assaying subspace $X$ will be equipped with the
Mackey topology $\tau(X,\overline{X^\#})$, unless it is otherwise stated.
The set of all continuous linear operators from $X$ into another
assaying subspace $Y$ will be denoted by $\cL(X, Y)$.

J.-P. Antoine and A. Grossmann \cite{Ant76a} introduced the operator spaces
$Op(V)$ (and more generally, $Op(V,W)$) of operators acting on $PIP$-spaces.
They also defined products for certain n-tuples of those operators.
In their paper J.-P. Antoine and A. Grossmann showed that,
given a non-degenerate $PIP$-space
$(V,\#,\langle . ,\! . \rangle)$, the space $Op(V)$ is
linearly isomorphic in a natural way
to the space $\cL(V^\#, V)$  of all continuous linear
operators mapping $V^\#$ into $V$. Using this isomorphism, we identify
here $Op(V)$ with $\cL(V^\#, V)$. Then the following definition of products
on $Op(V)$ is equivalent to the special case of the definition in
\cite{Ant76a}, where all operators to be multiplied belong to
$Op(V)=Op(V, V)$ for a fixed $PIP$-space $(V,\#,\langle.,\!.\rangle)$.
\begin{thd}             \label{circ}
Let $(V,\#,\langle.,.\rangle)$ be a non-degenerate $PIP$-space.
The factorization product $T_n\circ\dots\circ T_1$ of elements of
$Op(V)(=\cL(V^\#, V))$
is said to be defined if there are assaying subspaces $E_0, \dots ,
E_n $ of $V$ and continuous extensions $S_j \in {\cal L}(E_{j-1}, E_j)$
of $T_j$. In this case
$$T_n \circ \dots \circ T_1 \varphi= S_n(\dots(S_1 \varphi)\dots)\ .$$
\end{thd}

On $Op(V)$ there is defined an involution ${\rm A}\mapsto{\rm A}^\ast$
such that $\langle {\rm A}\,\varphi, \psi\rangle= \langle
\varphi,{\rm A}^\ast\,\psi\rangle$ for all $\varphi,\psi\in V^\#$,
i.e., ${\rm A}^\ast$ is the dual ${\rm A'}\in
\cL(\overline{V^\#},\overline{V})$ of A, considered as an element
of $\cL(V^\#,V)$.

 Our aim is to construct a linear mapping
$\pi \, : \, \A \rightarrow Op(V)$ that respects adjoints and
products, i.e., $\pi(x^\ast)=\pi(x)^\ast$ for all $x\in\A$ and
$\pi(xy)=\pi(x)\circ\pi(y)$ whenever $y\in R(\{x\})$.
As in \cite{AST95},
the basic idea is to define the operator $\pi(x)$ on the set
$\{b+\fN \,;\, b\in \B\}\subset V$ (see Proposition
\ref{PIPspace} for notations) by setting
$\pi(x)(b+\fN) = x b+\fN$.
Unfortunately, we are facing two difficulties:
In general, $\pi(x)$ is not yet defined on $V^\#$ since
$\{b+\fN \,;\, b\in \B\}$ is not necessarily equal to $V^\#$,
and for $y\in R(\{x\})$ we must find an assaying subspace
$X\subset V$ such that $\pi(y)\in \cL(V^\#, X)$ and that
$\pi(x)$ has an extension belonging to $\cL(X, V)$.
Under these circumstances it seems natural
that we need further assumptions. For this purpose we state the
following definition.
\begin{thd}             \label{hypostetig}
Suppose $\A$ is a partial $^\ast$-algebra,
$\B$ is a linear subspace of $R(\A)$, and $\Omega$ is
$\B$-weight on $\A$. If $X$ and $Y$  are linear subspaces of
$\B^\sharp$ and \mbox{$X\subset Y^\sharp$}, then $\Sigma(X,Y)$ denotes
the topology on $X$ that is generated by the family of semi-norms
$\{p_y\}_{y\in Y}$, where
$p_y(x)\stackrel{\rm def}{=}|\Omega(x,y)|$.
We say the partial product of $\A$
is $\Omega$-hypocontinuous w.r.t.\  $\B$, if the linear functionals
$$
\B\,\ni\, b\,\longmapsto\,\Omega(xb,w)\,\in\,\C
$$
are continuous for all $x\in \A$ and $w\in (x\B)^\sharp$ w.r.t.\
the topology $ \Sigma(\B,\B^\sharp)$, and if  for each  $x\in\A$
there exists a family $\cM$ of $\Sigma(\B^{\sharp\sharp},\B^\sharp)$-bounded
subsets of $\B$ such that $\B^{\sharp\sharp}=\cup_{M\in\cM}\overline{M}$,
where $\overline{M}$ denotes the closure of $M$ w.r.t.\
$\Sigma(\B^{\sharp\sharp},\B^\sharp)$, and the closed, convex, circled hull
of $xM$ is quasi-compact w.r.t.\
$\Sigma((x\B)^{\sharp\sharp},(x\B)^\sharp)$ for each $M\in\cM$.
\end{thd}
\noindent
{\bf Remarks: 1.} In general, the topology $\Sigma(X,Y)$ does not separate
the points of $X$. That's why we require the closed, convex, circled
hull of $xM$ to be quasi-compact; the term
``compact'' we reserve for compact Hausdorff spaces.\\
{\bf 2.} The terminology ``hypocontinuous'' alludes
to the concept of $\mathfrak{S}$-hypo\-con\-tin\-u\-ous bilinear forms as
defined by Bourbaki \cite[ch.III, $\S$5, 3.]{Bou87}.
How $\Omega$-hypo\-con\-tin\-u\-ous partial products
are related to $\mathfrak{S}$-hypocontinuous sesquilinear forms
will become clear in the proof of the next proposition.
\medskip

Now we are in a position to prove the following version of a generalized
$GNS$-representation.
%
%
\begin{thp}   \label{GNScirc}
Let $\Omega$ be a $\B$-weight on a partial $^\ast$-algebra $\A$.
Suppose that
${\fN}\stackrel{\rm def}{=} \{x\in\B^\sharp\,;\,\Omega(x,y)=0\
\mbox{ for all } y \in \B \}$ satisfies the assertions of
Proposition \ref{PIPspace}. Let $(V,\#,\langle.,\!.\rangle)$
denote the non-de\-gen\-e\-rate $PIP$-space defined in
Proposition \ref{PIPspace}. Assume that the partial product of $\A$
is $\Omega$-hypocontinuous w.r.t.\  $\B$.
Then there exists a unique linear mapping
$\pi \, : \, \A \rightarrow Op(V)$ such that
$\pi(x)(b+\fN) = x b+\fN$ for all $b\in\B$ and $x\in\A$.
Furthermore, $\pi$ satisfies
$\pi(x^\ast)=\pi(x)^\ast$ for all $x\in\A$ and
$\pi(xy)=\pi(x)\circ\pi(y)$ whenever $y\in R(\{x\})$.
\end{thp}

It is useful to introduce some temporary notations. Let
$\iota:\B^\sharp \rightarrow\B^\sharp/\fN$ denote the
canonical mapping. For an element
$\iota(a)\in\B^\sharp/\fN\,\mbox{(=V)}$
we write synonymously $a+\fN$ and $\hat{a}$; similarly we write
$\hat{M}$ in place of $\iota(M)$, where \mbox{$M\subset\B^\sharp$}.
The proof requires four steps: We shall define a linear mapping
\mbox{$\tilde{\pi}(x):\hat{\B}\rightarrow V$}, give a unique extension
$\pi(x)\in\cL(V^\#,V)$ of
$\tilde{\pi}(x)$, prove that \mbox{$\pi(x^\ast)=\pi(x)^\ast$},
and finally show that
$\pi(xy)=\pi(x)\circ\pi(y)$ whenever $y\in  R(\{x\})$.
\medskip\\
{\sc Proof:} {\it First step:}
As mentioned above, for each $x\in\A$
we define a linear mapping  $\tilde{\pi}(x):\hat{\B}\rightarrow V$
by setting
$$
\tilde{\pi}(x)(b+\fN)\,=\,xb+\fN\ .
$$
We have to show that this mapping is well defined. By
Definition \ref{Bweight} i), $xb\in \B^\sharp$ for all
$x\in\A$ and $b\in\B$. Given $b_1,b_2\in\B$ such that
$b_1-b_2\in\fN$, it follows from Definition \ref{Bweight} ii)
and $\fN=\fN_1$ that
$$ \Omega(xb_1-xb_2,b)\,=\,\Omega(b_1-b_2,x^\ast b)\,=\,0
\quad {\rm for\ all\ } b\in\B\ .
$$
This implies $xb_1-xb_2\in\fN$, hence $\tilde{\pi}(x)$ does not
depend on the choice of representatives. \\
{\it Second step:}
Our next goal is to show that $\tilde{\pi}(x)$ admits a unique
continuous extension $\pi(x):V^\#\rightarrow
(\stackrel{\wedge}{x\B})^{\#\#}$.
Let $\cS(V^\#,(\stackrel{\wedge}{x\B})^{\#})$ denote the set
of all separately continuous sesquilinear forms
$\beta:V^\#\times(\stackrel{\wedge}{x\B})^{\#}\rightarrow\C$\@.
The proof hinges on the (real linear) isomorphisms
$$
\cS(V^\#,(\stackrel{\wedge}{x\B})^{\#})\,\cong\,
   \cL(V^\#,(\stackrel{\wedge}{x\B})^{\#\#})\,\cong\,
\cL((\stackrel{\wedge}{x\B})^{\#},V)\ ,
$$
where $\cL(V^\#,(\stackrel{\wedge}{x\B})^{\#\#})\cong
\cL((\stackrel{\wedge}{x\B})^{\#},V)$
consists in
taking adjoints,
that is, given a $\beta\in\cS(V^\#,(\stackrel{\wedge}{x\B})^{\#})$,
there exist unique
${\rm B}\in \cL(V^\#,(\stackrel{\wedge}{x\B})^{\#\#})$ and
${\rm C}\in$ \mbox{$\cL((\stackrel{\wedge}{x\B})^{\#},V)$} such that
$$
\beta(\hat{v},\hat{w})\,=\,
\langle {\rm B}\,\hat{v}, \hat{w}\rangle\, =\,
\langle\hat{v},{\rm C}\, \hat{w}\rangle \quad
{\rm for\ all\ }\hat{v}\in V^\#{\rm \ and\ }
\hat{w}\in(\stackrel{\wedge}{x\B})^{\#}\ .
$$
The isomorphisms can be obtained by identifying the sesquilinear
forms on \mbox{$V^\#\times(\stackrel{\wedge}{x\B})^{\#}$}
with the bilinear forms on
$V^\#\times\overline{(\stackrel{\wedge}{x\B}})^{\#}$ and
applying the appropriate results for bilinear forms
(see K\"othe \cite[$\S$40, 1.]{Koe79}).

It is also known that a linear mapping is continuous w.r.t.\
the corresponding Mackey topologies if and only if
it is continuous w.r.t.\
the corresponding weak topologies. If we refer to the latter case,
we shall say ``weakly continuous''.

With the notations established in Definition \ref{hypostetig},
we observe that the weak topology
$\sigma(\hat{X},\overline{\hat{X}^\#})$ is the quotient topology of
$\Sigma(X,X^\#)$, where $X$ is an assaying subspace of $\B^\sharp$.
To see this, note that
\begin{eqnarray*}
\fN &=& \{x\in\B^{\sharp\sharp}\,;\,\Omega(x,y)=0\
\mbox{ for all } y \in \B^{\sharp} \}\\ &\subset& 
\{x\in X\,;\,\Omega(x,y)=0\ \mbox{ for all } y \in X^\sharp \}\\
&\subset& \{x\in\B^\sharp\,;\,\Omega(x,y)=0\
\mbox{ for all } y \in \B^{\sharp \sharp} \}\ =\ \fN\ ,
\end{eqnarray*}
hence $\fN=\{x\in X\,;\,\Omega(x,y)=0\ \mbox{ for all } y \in X^\sharp \}=
\cap_{y\in X^\sharp}\,p_y^{-1}(0)$. As a consequence,
$\iota:(X,\Sigma(X,X^\sharp))\rightarrow
(\hat{X},\sigma(\hat{X},\overline{\hat{X}^\#}))$ is a continuous,
open mapping. In addition, $(\hat{\B},\sigma(\hat{\B},\overline{V}))$ is
a topological subspace of $(V^\#,\sigma(V^\#,\overline{V}))$.

Next, for $x\in\A$ define
$$
\tilde{\beta}_x\,:\,\hat{\B}\times(\stackrel{\wedge}{x\B})^{\#}\,
\ni\, (\hat{b},\hat{w})\ \longmapsto\
\langle\tilde{\pi}(x)\,\hat{b},\hat{w}\rangle\,\in\,\C\ .
$$   Since
$\tilde{\pi}(x)\,\hat{b}\in(\stackrel{\wedge}{x\B})^{\#\#}$,
$\tilde{\beta}_x$ is weakly continuous in the second argument.
To prove continuity in the first argument, consider
the linear functionals
\begin{eqnarray*}
\Phi_{x,w}\,:\,\B\,\ni\,b &\longmapsto & \Omega(xb,w)\,\in\,\C\ ,
\quad {\rm and}\\
\phi_{x,w}\,:\,\hat{\B}\,\ni\,\hat{b} &\longmapsto &
\langle\tilde{\pi}(x)\,\hat{b},\hat{w}\rangle\,\in\,\C\ ,
\quad {\rm for}\ x\in\A{\rm \ and\ }w\in (x\B)^{\sharp}\ .
\end{eqnarray*}
Since $\langle\tilde{\pi}(x)\,\hat{b},\hat{w}\rangle=\Omega(xb,w)$,
we have $\Phi_{x,w}=\phi_{x,w}\circ\iota$ and
$\phi_{x,w}^{-1}(A)=\iota(\Phi_{x,w}^{-1}(A))$, where $A\subset\C$\@.
But $\Phi_{x,w}$ is continuous by the $\Omega$-hypocontinuity
of the partial product, and $\iota$ is open, hence $\phi_{x,w}^{-1}(U)$
is open for every open set $U\subset\C$\@. This implies the weak
continuity of $\phi_{x,w}$ and, moreover, the weak continuity of
$\tilde{\beta}_x$ in the first argument.
Let $\cM$ be a family of subsets of $\B$ satisfying the assumptions of
Definition \ref{hypostetig}. Set
$\hat{\cM}=\{\hat{M}\,;\,M\in\cM\}$.
As $\sigma(V^\#,\overline{V})$ is the quotient topology
of $\Sigma(\B^{\sharp\sharp},\B^\sharp)$, we have
$V^\#=\cup_{\hat{M}\in\hat{\cM}}\, \overline{\hat{M}}$, where
the closure $\overline{\hat{M}}$ of $\hat{M}$ is taken w.r.t.\
$\sigma(V^\#,\overline{V})$, and all $\hat{M}\in\hat{\cM}$
are bounded. Applying the facts that the closed, convex,
circled hull  of $xM$ is quasi-compact w.r.t.\
$\Sigma((x\B)^{\sharp\sharp},(x\B)^\sharp)$
(see Definition \ref{hypostetig}) and that $\iota$ is continuous,
one verifies readily that the closed, convex,
circled hull of $(\stackrel{\wedge}{xM})$
is compact w.r.t.\  \mbox{$\sigma((\stackrel{\wedge}{x\B})^{\#\#},
(\overline{\stackrel{\wedge}{x\B}})^\#)$}; hence the polar
$(\stackrel{\wedge}{xM})^\circ=
\{\hat{y}\in(\stackrel{\wedge}{x\B})^\#\,;\,
|\tilde{\beta}_x(\hat{m},\hat{y})|\leq 1\ \mbox{ for all }
\hat{m}\in\hat{M}\}$
is a 0-neighbourhood in $(\stackrel{\wedge}{x\B})^\#$.
Now, collecting the properties of $\tilde{\beta}_x$ and
$\hat{\cM}$, we observe that $\tilde{\beta}_x$ is
$\hat{\cM}$-hypocontinuous as defined by Bourbaki \cite{Bou87}.
This is the crucial observation. It follows by a theorem concerning
hypocontinuous bilinear mappings
(see \mbox{\cite[ch.III, $\S$5, 4.]{Bou87}})
that $\tilde{\beta}_x$ has a unique separately continuous extension
$$
\beta_x\,:\,V^\#\times(\stackrel{\wedge}{x\B})^\#\,
                                  \longrightarrow\,\C\ .$$
By the above mentioned isomorphisms, this implies that there
exist unique linear operators
${\rm X}\in\cL(V^\#,(\stackrel{\wedge}{x\B})^{\#\#})$ and
${\rm Y}\in\cL((\stackrel{\wedge}{x\B})^{\#},V)$ such that
\begin{equation}  \label{eq1}
\beta_x(\hat{v},\hat{w})\,=\,
\langle {\rm X}\,\hat{v}, \hat{w}\rangle\, =\,
\langle\hat{v},{\rm Y}\, \hat{w}\rangle \quad
{\rm for\ all\ }\hat{v}\in V^\#{\rm \ and\ }
\hat{w}\in(\stackrel{\wedge}{x\B})^{\#}\ .
\end{equation}
In particular,
$$
\langle \tilde{\pi}(x)\,\hat{b},\hat{w}\rangle \, =\,
\langle {\rm X}\,\hat{b}, \hat{w}\rangle \quad
{\rm for\ all\ }\hat{b}\in \hat{\B}{\rm \ and\ }
\hat{w}\in(\stackrel{\wedge}{x\B})^{\#}\ ,
$$
which gives ${\rm X}\lceil_{\!\hat{\B}}=\tilde{\pi}(x)$ since
$(\stackrel{\wedge}{x\B})^\#$ separates the points of
$(\stackrel{\wedge}{x\B})^{\#\#}$. Furthermore,
${\fN}=\{x\in\B^\sharp\,;\,\Omega(x,y)=0\ \mbox{ for all } y\in\B\}$
implies that the polar of $\hat{\B}$ taken in $V$ is $\{0\}$,
and this is equivalent to the density of $\hat{\B}$ in $V^\#$.
Thus the continuous extension of $\tilde{\pi}(x)$
is unique. As the embedding $(\stackrel{\wedge}{x\B})^{\#\#}
\hookrightarrow V$ is continuous, X can also be considered as
an element of $\cL(V^\#,V)$. Setting $\pi(x)\,\hat{v}=
{\rm X}\,\hat{v}$ for all $\hat{v}\in V^\#$,
we obtain the desired continuous extension
$\pi(x):V^\#\rightarrow V$ of $\tilde{\pi}(x)$. \\
{\it Third step:}
We observe that
$$
\langle \pi(x)^\ast\,\hat{b}_1,\hat{b}_2\rangle \, =\,
\langle \hat{b}_1,\pi(x)\,\hat{b}_2\rangle\, =\,
\langle \pi(x^\ast)\,\hat{b}_1,\hat{b}_2\rangle\quad
{\rm for\ all\ }\hat{b}_1,\hat{b}_2\in \hat{\B}
{\rm \ and\ }x\in\A\ .
$$
The first equality follows from the definition of the involution
on $Op(V)$; the second equality follows by using
$\pi(x)(b+\fN)=xb+\fN$ for all $b\in\B$, the def\-i\-ni\-tion of
$\langle . ,\! . \rangle$ and Definition \ref{Bweight} ii).
Since $\hat{\B}$ separates the points of $V$ and is dense in
$V^\#$, the above equation implies
$\pi(x)^\ast =\pi(x^\ast)$.\\
{\it Fourth step:}
Suppose $x,y\in\A$ and $y\in R(\{x\})$. An application of
Definition \mbox{\ref{Bweight} iii)} shows that
$(\stackrel{\wedge}{y\B})\subset(\stackrel{\wedge}{x^\ast\B})^\#$
which gives $(\stackrel{\wedge}{y\B})^{\#\#}
\subset(\stackrel{\wedge}{x^\ast\B})^\#$. Furthermore, the
embedding $(\stackrel{\wedge}{y\B})^{\#\#}\hookrightarrow
(\stackrel{\wedge}{x^\ast\B})^\#$ is continuous. Replacing
$x$ by $y$ in Equation (\ref{eq1}) and using the continuity of
the embedding $(\stackrel{\wedge}{y\B})^{\#\#}\hookrightarrow
(\stackrel{\wedge}{x^\ast\B})^\#$, we find an operator
${\rm S}_1\in\cL(V^\#,(\stackrel{\wedge}{x^\ast\B})^\#)$
such that ${\rm S}_1\,\hat{v}=\pi(y)\,\hat{v}$ for all
$\hat{v}\in V^\#$ (see the final part of the second step).
Replacing $x$ by $x^\ast$ in Equation (\ref{eq1}), we find an
operator ${\rm S}_2\in\cL((\stackrel{\wedge}{x^\ast\B})^\#,V)$
such that
\begin{equation}  \label{eq2}
\langle \pi(x^\ast)\,\hat{v}, \hat{w}\rangle\, =\,
\langle\hat{v},{\rm S}_2\, \hat{w}\rangle \quad
{\rm for\ all\ }\hat{v}\in V^\#{\rm \ and\ }
\hat{w}\in(\stackrel{\wedge}{x^\ast\B})^{\#}\ .
\end{equation}
Moreover, Equation (\ref{eq2}) implies
$\pi(x^\ast)^\ast={\rm S}_2\lceil_{\!V^\#}$ which gives, by
the third step, $\pi(x)\,\hat{v}={\rm S}_2\,\hat{v}$ for all
$\hat{v}\in V^\#$. Hence the factorization product
$\pi(x)\circ\pi(y)$ is defined and satisfies
$\pi(x)\circ\pi(y)\,\hat{v}={\rm S}_2({\rm S}_1\,\hat{v})$
for all $\hat{v}\in V^\#$.

Now, for all $\hat{b}_1,\hat{b}_2\in \hat{\B}$ we have
$$
\langle \pi(xy)\,\hat{b}_1,\hat{b}_2\rangle \, =\,
\langle \pi(y)\, \hat{b}_1,\pi(x^\ast)\,\hat{b}_2\rangle\, =\,
\langle {\rm S}_2({\rm S}_1\,\hat{b}_1),\hat{b}_2\rangle\, =\,
\langle \pi(x)\circ\pi(y)\,\hat{b}_1,\hat{b}_2\rangle\ ;
$$
the first identity is obtained by using
Definition \ref{Bweight} iii), and the second identity follows
from Equation (\ref{eq2}). Since $\hat{\B}$ separates the points
of $V$ and is dense in $V^\#$, this implies
$\pi(xy)=\pi(x)\circ\pi(y)$, and the proof is complete.
\medskip\\
{\bf Remarks: 1.}
Let $E$ and $F$ be locally convex spaces. A theorem by
Bourbaki \cite[ch.III, $\S$5, 3.]{Bou87} asserts that
if $F$ is barrelled, then every
separately continuous bilinear form from $E\times F$ into $\C$
is $\mathfrak{S}$-hypo\-con\-tin\-u\-ous for any family
$\mathfrak{S}$ of bounded subsets of $E$.
An examination of the proof of Proposition \ref{GNScirc}
shows that if the spaces $(\stackrel{\wedge}{x\B})^{\#}$ are
barrelled  for all $x\in\A$, then we can replace the hypothesis
that the partial product on $\A$ is $\Omega$-hypocontinuous
w.r.t.\ $\B$ by the statement that all linear functionals
$\B\,\ni\, b\,\longmapsto\,\Omega(xb,w)\,\in\,\C$
are continuous for all $x\in\A$ and \mbox{$w\in(x\B)^\sharp$}
w.r.t.\ $\Sigma(\B,\B^\sharp)$.\\
{\bf 2.} Let $\fN$ satisfy the assertions of
Proposition \ref{PIPspace}.
If $\B^{\sharp\sharp}=\B+\fN$, and if the linear functionals
$
\B\,\ni\, b\,\longmapsto\,\Omega(xb,w)\,\in\,\C
$
are continuous for all $x\in\A$ and \mbox{$w\in(x\B)^\sharp$}
w.r.t.\ $\Sigma(\B,\B^\sharp)$, then the partial product on $\A$
is at once $\Omega$-hypocontinuous w.r.t.\ $\B$.
To see this, set $\cM=\{b\}_{b\in\B}$ and note
that $\overline{\{b\}}=\{b+n\,;\,n\in\fN\}$.
\medskip

Using the preceding remarks, we can restate Proposition \ref{GNScirc}
in the following way.
\begin{thp}
Let $\Omega$ be a $\B$-weight on a partial $^\ast$-algebra $\A$ and
suppose that
${\fN}\stackrel{\rm def}{=} \{x\in\B^\sharp\,;\,\Omega(x,y)=0\
\mbox{ for all } y \in \B \}$ satisfies the assertions of
Proposition \ref{PIPspace}. 
Let $(V,\#,\langle.,\!.\rangle)$
denote the non-de\-gen\-e\-rate $PIP$-space defined in
Proposition \ref{PIPspace}.
Assume that the linear functionals
$\B\,\ni\, b\,\longmapsto\,\Omega(xb,w)\,\in\,\C$
are continuous for all $x\in\A$ and $w\in(x\B)^\sharp$
w.r.t.\  $\Sigma(\B,\B^\sharp)$.
Suppose that one of the following conditions is satisfied:\\
i) $\B^{\sharp\sharp}=\B+\fN$,\\
ii) $(\stackrel{\wedge}{x\B})^{\#}$ is barrelled for all $x\in\A$.\\
Then there exists a unique  linear mapping
$\pi$ from $\A$ into $Op(V)$ such that
$\pi(x)(b+\fN)=x b+\fN$ for all $b\in\B$ and $x\in\A$.
Furthermore, $\pi$ satisfies
$\pi(x^\ast)=\pi(x)^\ast$ for all $x\in\A$ and
$\pi(xy)=\pi(x)\circ\pi(y)$ whenever $y\in R(\{x\})$.
\end{thp}

Example \ref{ex2} will show that we cannot dispense with a hypothesis
that ensures in the preceding propositions the existence of
the products $\pi(x)\circ\pi(y)$. Here we required the partial product
of $\A$ to be $\Omega$-hypocontinuous w.r.t.\ $\B$. It should
be pointed out that this assumption is sufficient but we did not
prove that it is necessary; it seems to be rather difficult to give
a necessary condition. Observe that the factorization of the product
$\pi(x)\circ\pi(y)$ was achieved by proving that
$\pi(y)\in\cL(V^\#,(\stackrel{\wedge}{y\B})^{\#\#})$ and that
$\pi(x)$ admits an extension belonging to
$\cL((\stackrel{\wedge}{x^\ast\B})^{\#},V)$, but any assaying
subspace $X$ such that $(\stackrel{\wedge}{y\B})^{\#\#}
\subset X\subset(\stackrel{\wedge}{x^\ast\B})^{\#}$, that
$\pi(y)\in\cL(V^\#,X)$, and that $\pi(x)$ possesses an extension
in $\cL(X,V)$ factorizes the product $\pi(x)\circ\pi(y)$.
A necessary condition would have to control all those assaying
subspaces $X$.
%
%
\section{Representations based on extended products}
%
%

In this section, there will be  considerd
a second approach to $GNS$-rep\-re\-sen\-ta\-tions
of partial $^\ast$-algebras which uses a more general partial
product $T_1 \ast T_2$ of
operators on $PIP$-spaces. This product is similar to the weak product
defined in \cite{Ant85} for certain operators on Hilbert spaces. In contrary
to the product considered in Section 3, it is defined only for
certain pairs of operators, but not for n-tuples. However, it allows
to construct a representation of a partial $^\ast$-algebra $\A$ based on a
$\B$-weight $\Omega$, whenever besides the necessary conditions given
in Section 2 also the quite general additional condition
$\B^{\sharp \sharp}= \B+\fN$ is satisfied. An example shows that
it is impossible to define operators $\pi(x)\ \, (x \in \A)$ as
elements of $Op(V)$ in a natural way without any additional condition.
There will also be considered a further partial product $\bullet$, introduced 
in \cite{Wag97},  which
plays an intermediate role between the partial products $\circ$ and $\ast$.
Several examples demonstrate properties of these products.

The products $\bullet$ and $\ast$ are defined as follows.
%
%
\begin{thd}     \label{ast}
Let $(V, \#, \langle .,\!.\rangle)$ be a non-degenerate $PIP$-space.
The product $T_2\ast T_1$ of two elements of $Op(V)$ is defined if and only
if the following two equivalent conditions are satisfied:\\
i) There exists a $C\in Op(V)$ such that
$$
\langle T_1\varphi,T^\ast_2\psi\rangle\ =\ \langle C \varphi,\psi\rangle
\qquad(\mbox{for all } \varphi,\psi\in V^\#)\ .
$$
ii) There exist linear mappings $C,D:V^\#\rightarrow V$ such that
$$
\langle T_1\varphi,T^\ast_2 \psi\rangle\ =\ \langle C \varphi,\psi\rangle
\ =\ \langle \varphi, D \psi\rangle\qquad
(\mbox{for all } \varphi, \psi \in V^\#)\ .
$$
In this case
$$
T_2\ast T_1\,  \stackrel{\rm def}{=}\,  C\ .
$$
\end{thd}
%
%
\begin{thd}            \label{bullet}
Let $(V, \#, \langle.,\!.\rangle)$ be a non-degenerate $PIP$-space.
The product $T_2\bullet T_1$ of two elements of $Op(V)$
is defined if and only
if there exist assaying subspaces $X$, $Y$ of $V$ such that
the following four conditions are satisfied:\\
i) $T_1(V^\#) \subset X$,\\
ii) $T_2^\ast(V^\#) \subset Y$, \\
iii) $T_2$ has a continuous extension $S:X \rightarrow V$, \\
iv) $T_1^\ast$ has a continuous extension $R:Y \rightarrow V$.
In this case
$$
T_2 \bullet T_1 \varphi\ \stackrel{\rm def}{=}\ S(T_1 \varphi)\qquad
\mbox{for }
\varphi \in V^\#\,.
$$
\end{thd}

Note that $T_2 \bullet T_1$ belongs to $Op(V)$ since its adjoint is
given by
$$
(T_2 \bullet T_1)^\ast \varphi\ =\ R(T_2^\ast \varphi) \qquad
\mbox{for }
\varphi \in V^\#\,.
$$
Note also that the existence of $T_2\circ T_1$ implies
$T_2\bullet T_1=T_2\circ T_1$, and that the existence of $T_2\bullet T_1$
implies
$T_2\ast T_1=T_2\bullet T_1$. Actually, the domains of definition of
the three partial products are different. This may be seen by using
Example \ref{ex4}
below and Examples 3.5 and 4.4 in \cite{Kue98}.

Now it is not difficult to construct a representation by using the partial
product $\ast$
(cf. \cite{Kue98}).
The construction of the representation and some of its properties
are described in the following proposition.

%
%
\begin{thp}         \label{Prop4-3}
Let $\Omega$ be a $\B$-weight on a partial $^\ast$-algebra $\A$.
Suppose that $\fN$ is a linear subspace of $\B^\sharp$ such that
the assertions of Proposition \ref{PIPspace} are satisfied and that
$\B^{\sharp \sharp}=\B+\fN$. Let $(V, \#, \langle.,\!.\rangle)$ denote
the non-degenerate $PIP$-space defined in Proposition \ref{PIPspace} .
Then the formula
$$\pi(x)(b+\fN)\ =\ x b+\fN \qquad (b \in \B)$$
defines a linear mapping $\pi \, : \, \A \rightarrow Op(V)$ such that
$\pi(x^\ast)=\pi(x)^\ast$ for all $x \in \A$. Moreover
$$\pi(x_2 x_1)\ =\ \pi(x_2) \ast \pi(x_1)$$
for all $x_1, x_2 \in \A$ with $x_1 \in R(\{ x_2 \})$.
\end{thp}

Note that, provided assertion ii) in
Proposition \ref{PIPspace}  is satisfied,  condition
iv) in Definition \ref{Bweight}
means that $\{b+\fN \,;\, b\in \B\}$ is a dense linear
subspace of $\B^{\sharp\sharp}/\fN=V^\#$. By the following two examples,
this does not imply that there
exist operators $\pi(a)\in Op(V)$ satisfying $\pi(a)(b+\fN)=a b+\fN$ for all
$a\in \A$ and $b \in \B$.
%
%
\begin{thex}            \label{ex2}
There exist a partial $^\ast$-algebra $\A$, a linear subspace
$\B \subset \A$, and a $\B$-weight $\Omega$ on $\A$, such that the spaces
$\fN_1$ and $\fN_2$ defined in Proposition \ref{PIPspace} are equal to
$\{ 0 \}$ and
such that for some $a \in \A$ the mapping
$$
\B\, \ni\, x\, \rightarrow\, a \, x \,\in\, \B^\sharp \,(\,=V)
$$
does not extend to an element of $Op(V)$.
\end{thex}

Indeed, let $\A$ be the $^\ast$-algebra $\omega$ of all
complex valued sequences
(with pointwise algebraic operations), let $\B$ be the $^\ast$-subalgebra
$\varphi$ of all sequences of finite support, and let $\Omega$
be the usual scalar product of $l_2 \subset \omega$. Then, in particular,
$\Gamma(\sharp)=l_2 \times l_2$ and it is easy to see that $\Omega$
is a $\B$-weight. Moreover, using the notations of
Proposition \ref{PIPspace}, we have
$\fN_1=\fN_2=\{ 0 \}$ and $V \,(=\B^\sharp)=
V^\# = l_2$. Consequently, $Op(V)$ coincides with the space of
all bounded linear operators on $l_2$ and the example is completed by
taking $a=(a_j)_{j \in \N}$ to be an unbounded sequence.

%
%
\begin{thex}               \label{ex3}
There exist a partial $^\ast$-algebra $\A$, a linear subspace
$\B \subset \A$, and a $\B$-weight $\Omega$ on $\A$, such that
the equivalent conditions
of Proposition \ref{PIPspace} are satisfied and
that for some $a \in \A$  the space
$\fN_1 = \fN_2$ defined in Proposition \ref{PIPspace} is
not invariant for the mapping
$$
\B\, \ni\, x\, \rightarrow\, a \, x\, \in \B^\sharp\ .
$$
\end{thex}

This example will be constructed in the space
$$
\A\ \stackrel{\rm def}{=}\ \C^5\ =\ \{ (x_1, x_2, x_3, x_4, x_5);\,
x_j \in \C \}\ ,
$$
endowed with its usual structure of a $^\ast$-vector space and with the
(commutative but non-associative) partial product defined as follows.
\begin{eqnarray*}
\fC & \stackrel{\rm def}{=} & \{ (x_j)_{j=1}^5 \in \A ;\, x_5=0 \}\ , \\
\Gamma  & \stackrel{\rm def}{=} & (\fC \times \A) \cup (\A \times \fC)\ ,\\
(x_j)_{j=1}^5 \cdot (y_j)_{j=1}^5 &  \stackrel{\rm def}{=} &
(x_1 y_1, x_2 y_2, x_3 y_3, x_4 y_4, 0)+ (x_5 y_3, 0, x_5 y_1, 0, 0) \\  & &
+(x_3 y_5, 0, x_1 y_5, 0, 0) \\ & & \mbox{for }
((x_j)_{j=1}^5, (y_j)_{j=1}^5) \in \Gamma\ .
\end{eqnarray*}
Furthermore, we define a $\B$-weight by setting
\begin{eqnarray*}
\B & = &  \{ (x_j)_{j=1}^5 \in \A ;\,
x_3=x_4=x_5=0 \}\ , \\
\Gamma(\sharp) & = & \fC \times \fC\ , \\
\Omega((x_j)_{j=1}^5, (y_j)_{j=1}^5) & = & \sum_{j=2}^4 \, x_j \,
\overline{y_j}
\qquad \mbox{ for }
((x_j)_{j=1}^5, (y_j)_{j=1}^5) \in \Gamma(\sharp)\ .
\end{eqnarray*}
Easy computations show that all properties of
Definition \ref{Bweight} are satisfied.
E.g., it follows for $x=(x_j)_{j=1}^5 \in \A$,
$y=(y_j)_{j=1}^5 \in \fC$,
and $b=(b_j)_{j=1}^5, c=(c_j)_{j=1}^5 \in \B$ that
$y c \in \B$ and that
$$
\Omega(xb, yc)\ =\  x_2 \,b_2\,\overline{y_2}\,\overline{c_2}\ =\
\Omega(b, (x^\ast y)c)\ .
$$
By commutativity, this implies Definition \ref{Bweight} iii).

Clearly $\B^\sharp=\B^{\sharp \sharp} =\fC$ and
$\fN_1 = \fN_2= \C \cdot (1, 0, 0, 0, 0)$. Setting $a=(0, 0, 0, 0, 1)$,
we have $a \cdot (1, 0, 0, 0, 0) \not\in \fN_1$, which completes the example.
\medskip

%
%
Our final example shows that the representation described in
Proposition \ref{Prop4-3}
cannot be constructed by using the product on $Op(V)$ defined in
Definition \ref{circ}, in general.
\begin{thex}                \label{ex4}
There exist a partial $^\ast$-algebra $\A$, a linear subspace
$\B\subset\A$, and a $\B$-weight $\Omega$ on $\A$, such that all assumptions
of Proposition \ref{Prop4-3} are satisfied and that for some $a \in \A$ with
$a \in R(\{ a \})$ the product $\pi(a) \circ \pi(a)$ does not exist in the
sense of Definition \ref{circ}.
\end{thex}

Let again $\omega$ denote the $^\ast$-algebra of all complex valued sequences
$(x_j)=(x_j)_{j \in \N}$
(with pointwise operations). Let $\varphi$ be the $^\ast$-subalgebra
of all sequences of finite support. Consider elements $a = (a_n)$
and $a^2=((a_n)^2)$ of $\omega$, where $(a_n)_{n \in \N}$ is
a fixed unbounded sequence of positive real numbers. Define now
\begin{eqnarray*}
\A &=& \varphi + \C a + \C a^2\  , \\
\Gamma &=& (\varphi+\C a) \times (\varphi+\C a) \cup \A \times \varphi \cup
        \varphi \times \A\  .
\end{eqnarray*}
Endowed with the linear operations and the involution induced from $\omega$
and with the partial product obtained as the restriction of the product of
$\omega$ to $\Gamma$, $\A$ is a partial $^\ast$-algebra. Setting
$\B = \varphi$ and defining $\sharp$ and $\Omega$ by
\begin{eqnarray*}
\Gamma(\sharp) &=& \A \times \varphi \cup \varphi \times \A \ , \\
\Omega((x_j) , (y_j)) &=& \sum_{j=1}^\infty x_j \overline{y_j} \qquad
        ((x_j) , (y_j) \in \Gamma(\sharp))\ ,
\end{eqnarray*}
we get a $\B$-weight $\Omega$ on $\A$. It is easy to see that the
assertions of Proposition \ref{PIPspace} are satisfied for $\fN=\{ 0 \}$ and that
$\B^{\sharp \sharp}=\B$. Consequently, the $PIP$-space
$(V, \# , \langle.,\!.\rangle)$ constructed in Proposition \ref{PIPspace}
coincides here with $(\A , \sharp , \Omega( . , . ))$, and all assumptions
of Proposition \ref{Prop4-3} are satisfied. Given $(x_j) \in \A$, the operator
$\pi((x_j))$ acts on $V^\#=\B $ as multiplication operator with the
sequence $(x_j)$.

We show that $\pi(a) \, : \, V^\# \rightarrow V^\#$ is not continuous.
Otherwise it would have a continuous adjoint $A \,: \, V \rightarrow V$ such
that
$$\langle \pi(a) b_1 , b_2 \rangle\ =\ \langle b_1 , A \, b_2 \rangle\qquad
        (\mbox{for all } b_1 \in V^\# \mbox{ and } b_2 \in V)\ .$$
This would imply that $A \, a^2 = ((a_j)^3)$, which is impossible.

It can be seen in the same way that $\pi(a)$ does not have a continuous
extension to an operator $A \,: \, V \rightarrow V$. Since the only
assaying subspaces are $V$ and $V^\#$, this implies that the product
$\pi(a) \circ \pi(a)$ does not exist in the sense of Definition \ref{circ}.
Clearly, $\pi(a) \ast \pi(a) = \pi(a^2)$ by Proposition \ref{Prop4-3}.
\medskip\\
{\bf Remarks: 1.}
Since in the previous example all operators $\pi(x)$
($x \in \A$) satisfy $\pi(x) V^\# \subset V^\#$, it can be seen
easily that in this example  the equation
$$\pi(x_2 x_1)\ =\ \pi(x_2) \bullet \pi(x_1)$$
is satisfied for all $x_1, x_2 \in \A$ with $x_1 \in R(\{ x_2 \})$.
However, such a property is not satisfied in the general case, as it can
be shown by using Example 4.4 in  \cite{Kue98}.
\\
{\bf 2.} In particular it follows from Example \ref{ex4} that there exist
operators $T_1$ and $T_2$ on some $PIP$-space such that $T_2 \bullet T_1$
exists, but $T_2 \circ T_1$ is not defined. Similarly, it can be shown by
using Examples 3.5 or 4.4 in \cite{Kue98}
that there exist
operators $T_1$ and $T_2$ on some $PIP$-space such that $T_2 \ast T_1$
exists, whereas $T_2 \bullet T_1$ is not defined.


\end{document}